\documentclass[fleqn]{mat01}
\usepackage{times,mathtimy,amssymb,latexsym,amscd}
\begin{document}

\setcounter{page}{429}
\firstpage{429}

\newcommand{\A}{\mathbb{A}}
\newcommand{\R}{\mathbb{R}}
\newcommand{\reals}{\mathbb{R}}
\newcommand{\C}{\mathbb{C}}
\newcommand{\N}{\mathbb{N}}
\newcommand{\D}{\mathbb{D}}
\newcommand{\Z}{\mathbb{Z}}
\newcommand{\integers}{\mathbb{Z}}
\newcommand{\Q}{\mathbb{Q}}
\newcommand{\HH}{\mathbb{H}}
\newcommand{\K}{\mathbb{K}}
\newcommand{\E}{\mathbb{E}}
\newcommand{\F}{\mathbb{F}}

\newtheorem{theore}{Theorem}
\renewcommand\thetheore{\arabic{section}.\arabic{theore}}
\newtheorem{theor}[theore]{\bf Theorem}
\newtheorem{lem}[theore]{\it Lemma}
\newtheorem{propo}[theore]{\rm PROPOSITION}
\newtheorem{coro}[theore]{\rm COROLLARY}
\newtheorem{definit}[theore]{\rm DEFINITION}
\newtheorem{probl}[theore]{\it Problem}
\newtheorem{exampl}[theore]{\it Example}

\title{Transversals in non-discrete groups}

\markboth{Ramji Lal and R~P Shukla}{Transversals in non-discrete
groups}

\author{RAMJI LAL and R~P~SHUKLA}

\address{Department of Mathematics, University of Allahabad,
Allahabad~211~002, India\\
\noindent E-mail: ramjilal@mri.ernet.in; rps@mri.ernet.in}

\volume{115}

\mon{November}

\parts{4}

\pubyear{2005}

\Date{MS received 2 August 2004; revised 4 August 2005}

\begin{abstract}
The concept of `topological right transversal' is introduced to
study right transversals in topological groups. Given any right
quasigroup $S$ with a Tychonoff topology $T$, it is proved that
there exists a Hausdorff topological group  in which $S$ can be
embedded algebraically and topologically  as a right transversal
of a subgroup (not necessarily closed). It is also proved that if
a topological right transversal $(S, T_{S}, T^{S}, \circ)$ is such
that $T_{S} = T^{S}$ is a locally compact Hausdorff topology on
$S$, then $S$ can be embedded as a right transversal of a closed
subgroup in a Hausdorff topological group which is universal in
some sense.
\end{abstract}


\maketitle

\section{Introduction}

The importance of the study of transversals in groups individually
and abstractly has been explained in \cite{ra:tran},
\cite{sh:perf}, and \cite{shu:con}. Every right transversal of a
subgroup in a group has the structure of a right quasigroup with
identity. Conversely, every right quasigroup with identity can be
embedded as a right transversal into a group which is universal in
some sense (see Theorem~3.4 of \cite{ra:tran}). Transversals in a
group have more control on the structure of the group than
subgroups. For example, there are uncountably many non-Dedekind
(2-generator infinite simple groups (Theorem~28.7 of
\cite{ols:gps}) and (Theorem~28.8 of \cite{ols:gps})) groups all
of whose non-trivial proper subgroups are isomorphic, whereas
there is no non-Dedekind group with all non-trivial proper
transversals isomorphic (as right quasigroups)
\cite{sh:perf,ra:tar}.

The abstract description of right transversals of a closed
subgroup in a topological group is a difficult problem and needs
more serious thought. In this paper, we make an attempt in this
direction.

\section{Right transversals in topological groups}

Let $S$ be a right transversal of a closed subgroup $H$ in a
Hausdorff topological group $G$. Algebraically, $S$ is a right
quasigroup with identity (the right quasigroup structure $\circ$
on $S$ is defined by $\{ x \circ y \} = Hxy \cap S$). Let $\langle
S \rangle$ denote the subgroup of $G$ generated by $S$ and $H_S$
denote the subgroup $\langle S \rangle\cap H$. Then $H_S = \langle
\{xy(x\circ y)^{-1}\hbox{:}\ x, y \in S\}\rangle$, and $H_SS =
\langle S \rangle$, where $\circ$ is defined as above. Identifying
$S$ with the set $H\backslash G$ of right cosets of $H$ in $G$, we
obtain a permutation representation $\phi\hbox{:}\
G\rightarrow \hbox{Sym}(S)$ ($\{\phi(g)(x)\}= S\cap Hxg, g\in
G$, $x\in S$ and we adopt the convention that $(r.s)(x) =
s(r(x))$, $r,s\in \hbox{Sym}(S) (x\in S$). Let $G_S =
\phi(H_S)$. Then $G_S$ depends only on the right quasigroup
structure $\circ$ on $S$ and not on the subgroup $H$ (see
\cite{ra:tran}). This group is called the group torsion of the
right quasigroup $(S, \circ)$ \cite{ra:tran}. It is easy to
observe that $\phi$ is injective on $S$ and if we identify $S$
with $\phi(S)$, then $\phi(\langle S \rangle) = G_SS$ is a group
which also depends only on the right quasigroup $(S, \circ)$ and
in which $S$ is a right transversal of $G_S$. This group has some
universal property \cite{ra:tran} and is called the group
extension of $(S, \circ)$. Also $\langle S \rangle = H_SS$ is a
topological subgroup of $G$ and $G_SS$ has quotient topological
group structure in which $G_S$ is closed (for $H_S =
\phi^{-1}(G_S) = H\cap \langle S \rangle$ is closed in $\langle S
\rangle$) and $S$ is a right transversal of $G_S$ in $G_SS$. We
have two Tychonoff topologies on $S$, viz. $T_{S}$, the subspace
topology on $S$ induced by the topology of $G_SS$ and $T^{S}$, the
quotient topology on $S$ induced by the natural projection from
$G_SS$ to $S$. Clearly $T_S$ is finer than $T^S$. Let $T_{G_S}$
denote the subspace topology on $G_S$ (with respect to which it is
a topological group) and $T^{G_S}$, the quotient topology on $G_S$
induced by the natural projection from $G_SS$ to $G_S$. The
topology $T^{G_S}$ need not be even a $T_1$ topology
(Example~2.3). Since $(S, \circ)$ is a right quasigroup with
identity, we have another map $\chi\hbox{:}\ S\times
S\rightarrow S$\break given by
\begin{equation*}
\chi(x,y)\circ  x = y.
\end{equation*}
The binary operations $\circ$ (and also $\chi$) may not be
continuous with respect to either of the topologies $T_S$ and
$T^S$ on $S$. For example take $G = SO(n)$ and $H = SO(n-1)$ where
$n\not\in\{1,2,4,8\}$. The quotient space $H\backslash G$ is
homeomorphic to $S^{n-1}$. There is a right transversal $S$ of $H$
in $G$ (see Example~2.3) such that $G_S = H$ and $G_SS = G$. Then
$(S,T^S)$ is homeomorphic to $S^{n-1}$ with the usual topology.
Since there is no continuous right (left) quasigroup structure on
$S^{n-1}$ for $n\not \in \{1,2,4,8\}$ \cite{ad:hinv}, $\circ$ and
$\chi$ are not continuous from $(S \times S, T^{S} \times T^{S})$
to $(S, T^{S})$, where $(S \times S, T^{S} \times T^{S})$ is the
product of $(S, T^{S})$ with itself. Next, there is a choice of a
right transversal of $S$ in $G$ (see Example~2.3 with ${\bf H}
=\R^n$) for which $T_S$ is non-homogeneous, and as such
$\circ\hbox{:}\ (S \times S, T_{S} \times T_{S})\rightarrow
(S,T_S)$ is not continuous. However, $\circ$ and $\chi$ are
continuous from $(S \times S, T_{S} \times T_{S})$ to $(S,
T^{S})$.

The above discussion prompts us to have the following definition.

\begin{definit}\label{2.1}$\left.\right.$\vspace{.5pc}

\noindent{\rm A quadruple $(S, T_{S}, T^{S}, \circ)$, where $(S,
\circ)$ is a right quasigroup with identity, $T_{S}$ and $T^{S}$
are two Tychonoff topologies with $T_{S} \supseteq T^{S}$ and
$T^S$ homogeneous is called a {\it topological right transversal}
if
\begin{enumerate}
\renewcommand\labelenumi{(\roman{enumi})}
\leftskip .15pc
\item $\circ$ is  continuous map from $(S \times S, T_{S} \times
T_{S})$ to $(S, T^{S})$, and also

\item the map $\chi\hbox{:}\ (S \times S, T_{S} \times T_{S})
\rightarrow (S, T^{S})$ defined by the equation $\chi(x,
y)\circ x = y$, is continuous.
\end{enumerate}}
\end{definit}

We shall say that a topological right transversal $(S, T_{S},
T^{S}, \circ)$ is   proper if $(S, \circ)$ is not a group. This
terminology is justified (for such a right transversal in any
group is always a proper transversal). Every choice of a right
transversal of a closed subgroup $H$ of a topological group $G$
determines a topological right transversal. Trivially, an improper
topological transversal $(G, T_{G}, T^{G}, \circ)$, where $T_G$ is
a discrete topology on $G$ and $T^G$ is any non-discrete
homogeneous Tychonoff topology cannot arise as a right
transversal in the above manner. Indeed, it comes from a
topological right transversal if and only if $T_G = T^G$.

The following example says that even a proper topological right
transversal need not arise as a right transversal of a closed
subgroup in a topological group.

\begin{exampl}{\rm \label{2.2}
For each sufficiently large odd $n$ ($n>10^{10}$) and $m>1$, there
is a central extension $K$ of $B(m,n)$ (the Burnside group of
exponent $n$ on a set containing $m$ elements) by a cyclic group
$C_n$ of order $n$ which admits no non-discrete topology
(Theorem~31.5 of \cite{ols:gps}). Take  $G = K \times D$, where
$D$ is a finite group having a non-normal subgroup $H$ of prime
order $p$. Let $T$ be a topology on $G$ with respect to which it
is a $T_1$-topological group. As observed earlier, the induced
subspace topology on $K\times \{e\}$ is the discrete topology
(where $e$ is the identity element in $D$). Since every discrete
subgroup of a $T_1$ topological group is closed (Theorem~5.10 of
\cite{h-r:har}), $K\times \{e\}$ is a closed subgroup of $G$. But,
then $K\times \{d\}$ is closed for each $d\in D$. Since $D$ is
finite, this implies $K\times \{d\}$ is open for each $d\in D$.
Since $K\times \{d\}$ is open and discrete for all $d\in D$, it
follows that $T$ is the discrete topology. This shows that $G$
also does not admit any non-discrete $T_1$-topology. Let $H$ be a
non-normal subgroup of $D$ of order $p$. Then there is a right
transversal $S$ of $\{e\}\times H$ in $G$ which generates $G$ and
$G_SS \cong G$ (see Proposition~3.10 of \cite{ra:tran}). Thus
$G_SS$ admits only discrete topology. Hence  $(S, T_{S}, T^{S},
\circ)$, where $\circ$ is the binary operation induced on $S$ with
respect to which it is a right quasigroup, and arises from a
topological group structure on $G_SS$ if and only if $T_S = T^S$
is the discrete topology. Take any non-discrete homogeneous
Tychonoff topology $T^S$ on $S$ (note that such a topology exists,
viz. the topology induced by the topology of $Q$ (by fixing a
bijection between $S$ and $Q$)), then $(S, T_{S}, T^{S}, \circ)$,
where $T_S$ is the discrete topology on $S$, is a topological
right transversal. Observe that $(S, T_{S}, T^{S}, \circ)$ will
not arise from any topological group structure on $G_SS$ (note
that $(S, T_{S}, T^{S}, \circ)$ is a proper transversal).}
\end{exampl}

\begin{exampl}\label{2.3}
\hskip -.3pc{\it Transversal structures on spheres in Hilbert
spaces.} \ \ {\rm Let ${\bf H}$ be a real Hilbert space with
$\langle\, , \rangle$ as inner product structure. Let $O({\bf H})$
denote the topological group of isometries of ${\bf H}$. For the
product in $O({\bf H})$, we adopt the convention that $(T_1T_2)(x)
= T_2(T_1(x))$; $T_1, T_2\in O({\bf H})$ and $x\in {\bf H}$. Let
$S^{\bf H} = \{x\in {\bf H}\hbox{:}\ \|x\| = 1\}$ denote the unit
sphere in ${\bf H}$ (in the case of ${\bf H} = \R^n$, we denote it
by $S^{n-1}$). Fix a point $e_0\in S^{\bf H}$. Let $V = \langle
e_0 \rangle^{\perp}$ (the orthogonal complement of the subspace
$\langle e_0 \rangle$ generated by $e_0$) and $O(V)= \{T\in O({\bf
H})\hbox{:}\ T(e_0) = e_0\}$. Then $O(V)$ is a closed subgroup of
$O({\bf H})$. The topological group $O({\bf H})$ acts transitively
on $S^{\bf H}$ from right and the isotropy group at $e_0$ is
$O({\bf H})_{e_0} = O(V)$. This gives us a bijective continuous
map $\phi$ from $O(V)\backslash O({\bf H})$, the quotient space of
right cosets of $O(V)$ in $O({\bf H})$ to $S^{\bf H}$ given by
$\phi(O(V)T) = T(e_0)$. It is easy to observe that $\phi$ is a
homeomorphism.

Let $x\in S^{\bf H}$. Define a map $R_{x}\hbox{:}\ {\bf
H}\rightarrow {\bf H}$ by
\begin{equation*}
R_x = \begin{cases} J\left(2 P_{\frac{e_0+x}{\|e_0+x\|}} -I
\right),
&\hbox{if}\: x \not = -e_0\\
-I, &\hbox{if}\: x = -e_0
\end{cases}
\end{equation*}
where for $a\in {\bf H}$, $P_a$ denotes the projection of ${\bf H}$
onto $\R a$ and $J$ is defined by
\begin{equation*}
J(x) = \langle x,e_0 \rangle e_0-(x- \langle x,e_0 \rangle e_0) =
2 \langle x,e_0 \rangle e_0-x.
\end{equation*}
Thus for any non-zero $a\in {\bf H}$,
\begin{equation*}
P_a(x) = \left\langle x,{\frac{a}{\|a\|}} \right\rangle
{\frac{a}{\|a\|}}, \quad x\in {\bf H}.
\end{equation*}
Clearly the map $R_x\in O({\bf H})$. Also the map $R\hbox{:}\
S^{\bf H}\rightarrow O({\bf H})$ defined by $R(x) = R_x$ is
injective and is a set theoretic section (not continuous at $-e_0$)
of the quotient map $\phi$. Thus the image $R(S^{{\bf H}}) =
\{R_x\hbox{:}\ x\in S^{\bf H}\}$ is a right transversal of $O(V) =
O({\bf H})_{e_0}$ in $O({\bf H})$. Identifying $R(S^{\bf H})$ with
$S^{\bf H}$ through the map $R$, we get a subspace topology
$T_{S^{\bf H}}$ on $S^{\bf H}$. We observe that the usual subspace
topology of $S^{{\bf H}}$ in ${\bf H}$ is same as the quotient
topology $T^{S^{\bf H}}$ on $S^{\bf H}$ induced by the map $\phi$.
This gives us a topological right transversal structure $(S^{\bf
H}, T_{S^{\bf H}}, T^{S^{\bf H}}, \circ)$ on $S^{\bf H}$. The
right quasigroup structure $\circ$ on $S^{\bf H}$ is defined by
\begin{equation*}
x\circ y = \begin{cases} \left(2P_{\frac{e_0+y}{\|e_0+y\|}}
-I\right)(J(x)), &\hbox{if}\: y \not = -e_0 \\
-x, &\hbox{if}\: y = -e_0
\end{cases}.
\end{equation*}
Since the section $R$ is not continuous at $-e_0$, the topology
$T_{S^{\bf H}}$ is strictly finer than $T^{S^{\bf H}}$ and is not
homogeneous. In fact $\{-e_0\}$ is an isolated point of $(S^{\bf
H}, T_{S^{\bf H}})$ and the subspace topologies on $S^{\bf
H}\!\setminus\!\{-e_0\}$ induced by $T_{S^{\bf H}}$ and $T^{S^{\bf
H}}$ are same. Thus, $(S^{\bf H}, T_{S^{\bf H}})$ is locally
compact. It can also be seen that $O({\bf H})$ is generated by
$S^{\bf H}$ (more precisely $R(S^{\bf H})$) and the core of $O(V)$
in $O({\bf H})$ is trivial. This implies that (see
Proposition~3.10 of \cite{ra:tran}) $G_{\!S^{\bf H}}S^{\bf H} =
O({\bf H})$. We also observe that the quotient topology
$T^{G_{\!S^{\bf H}}}$ induced by the projection map from $G_{\!S^{\bf
H}}S^{\bf H}$ to $G_{\!S^{\bf H}}$ is not a $T_1$-topology (for
$S^{\bf H}$ is not closed in $O({\bf H})$).}
\end{exampl}

\section{Transversals on locally compact Hausdorff
spaces}\vspace{.2pc}

\setcounter{theore}{0}

\begin{definit}{\rm \cite{ar:hom}.}\label{3.1}$\left.\right.$\vspace{.5pc}

\noindent {\rm Let $X$ be a locally compact Hausdorff space. Let
$G = {\rm Homeo}(X)$ be the group of homeomorphisms of $X$. Let
$K$ and $W$ be closed and open subspaces of $X$ respectively such
that either $K$ is compact or the complement of $W$ in $X$ is
compact. Let $(K,W)$ denote the set of all homeomorphisms of $X$
which send $K$ into $W$. The $g$-topology on $G$ is the topology
generated by the set $\cal S$ as a subbase whose members are of
the form $(K,W)$.}
\end{definit}

\begin{theor}[\!]\label{3.2}
Let $(S,T,T,\circ)$ be a topological right transversal{\rm ,}
where $T$ is a locally compact Hausdorff topology on $S$. Then
there exists a unique Hausdorff topological group $G^S$ containing
$(S,T,T,\circ)$ as a  topological right transversal such that
given any topological group $G$  containing $(S,T,T,\circ)$ as a
right transversal{\rm ,} there exists a unique continuous
homomorphism from $G$ to $G^S$ which is the identity map on $S$.
\end{theor}

\begin{proof}
Let $(S,T,T,\circ)$ be a topological right transversal with $T$ a
locally compact Hausdorff topology on $S$. Then the binary
operation $\circ$ and the map $\chi$ defined by $\chi(x,y)\circ x
= y$ are continuous from $(S\times S, T\times T)$ to $(S,T)$.
$S\!\setminus\!\{e\}$ ($e$ is the identity of $(S,\circ)$) being open
subspace of $S$ is also locally compact. Thus
Homeo $(S\!\setminus\!\{e\})$ is a topological group with respect to
the $g$-topology (Theorem~3 of \cite{ar:hom}). Also the natural
action $\theta\hbox{:}\ S \times$ Homeo$(S) \rightarrow S$ is
continuous \cite{di:trans}. Hence $\theta \vert S \times H^s =
\theta^{s}$, from $S \times H^s$ to $S$ is also continuous, where
$H^s = \{ f \in {\rm Homeo}(S)\hbox{:}\  f(e) = e \}$. Then $(S, H^s, \sigma^{s},
f^{s})$ is a $c$-groupoid (see Definition~2.1 and p.~78 of
\cite{ra:tran}), where $f^s$ is the map from $S\times S$ to $H^s$
defined by
\begin{equation*}
f^s(x,y)(z) = \chi(x\circ y, (z\circ x)\circ y), \quad x,y,z\in S
\end{equation*}
and the map $\sigma^s$ from $S\times H^s$ to $H^s$ defined by
($\sigma^s(x,h)$ is denoted by $\sigma_{x}^s(h)$),
\begin{equation*}
\sigma_{x}^s(h)(y) = \chi(x\theta^sh, (y\circ x)\theta^sh),\quad
x,y\in S, \quad h\in H^s.
\end{equation*}
Let $G^S = H^sS$ be the general extension of $H^s$ by $S$
determined by the $c$-groupoid $(S, H^s, \sigma^{s},f^{s})$ (p.~79
of \cite{ra:tran}). The product in $G^S$ is given by
\begin{equation*}
hx \cdot ky = h\sigma_{x}^s(k)f^s(x\theta^sk, y)x\theta^sk\circ y,
\end{equation*}
and the inverse of $hx$ is given by
\begin{equation*}
(hx)^{-1} = (f^s(x',x))^{-1}\sigma_{x'}^s(h^{-1})x'\theta^s
h^{-1},
\end{equation*}
where $h,k\in H^s$, $x,y\in S$ and $x'$ is the left inverse of $x$
in $(S, \circ)$ (p.~72 of \cite{ra:tran}). We note that the
subgroup $\langle S \rangle $ of $G^S$ generated by $S$ is $G_SS$.

Now, we show that $G^S$ is a topological group with respect to the
product topology on $H^s \times S$. The evaluation map
$E(f^{s})\hbox{:}\ (S \times S) \times S \rightarrow S$
(defined by $E(f^{s})((x, y), z) = f^{s}(x, y)(z)$), being the
compsosition of continuous maps $((x, y), z) \mapsto (x \circ y,
(z \circ x) \circ y)$ from $(S \times S) \times S$ to $S \times S$
and $\chi$, is continuous. It follows (Theorem~3.1, p.~261 of
\cite{du:topo}) that $f^{s}\hbox{:}\ S \times S \rightarrow
H^s$ is a continuous map. Consider the map
\begin{equation*}
\Sigma\hbox{:}\ S \times H^s \rightarrow H^s,  (x,h) \mapsto
\sigma^{s}_{x}(h).
\end{equation*}
Again the evaluation map
\begin{equation*}
E(\Sigma)\hbox{:}\ (S \times H^s) \times S \rightarrow S,
((x,h), y) \mapsto y \theta^{s} \sigma^{s}_{x}(h),
\end{equation*}
being the composition of the continuous maps $((x, h), y) \mapsto
(x \theta^{s} h, (y \circ x) \theta^{s} h)$ from $(S \times H^s)
\times S$ to $S \times S$ and $\chi$, is continuous. It follows
(Theorem~3.1, p.~261 of \cite{du:topo}) that $\Sigma$ is also
continuous. Further, the map $x \mapsto x'$ ($x'$ is the left
inverse of $x$ in $S$) from $S$ to $S$ being the composition of
continuous maps $x \mapsto (x, e)$ from $S$ to $S \times \{e \}$
and $\chi$, is continuous. Using the continuity of the above maps,
we find that the binary operation in $G^S$ and the map $g \mapsto
g^{-1}$ from $G^S$ to itself are continuous. Finally, we observe
that induced subspace topology on $S$ in $G^S$ is the given
topology on $S$ and induced right quasigroup structure on the
right transversal $S$ of the closed subgroup $H^s$ in $G^S$ is the
given right quasigroup structure.

Further, suppose that $G$ is a topological group such that
$(S,T,T,\circ)$ is a topological right transversal of a closed
subgroup $H$ in $G$. Then the action $\theta$ of $H$ on $S$
defined by $\{\theta(x,h)\} = S\cap Hxh$, where $x\in S$ and $h\in
H$, is continuous. As such it induces a continuous homomorphism
$q\hbox{:}\ H\rightarrow H^s$. Let $\phi\hbox{:}\
G\rightarrow G^S$ be the map defined by $\phi (ax) = q(a)x$.
Then $\phi$ is a group homomorphism (see Theorems~2.6 and 3.4 of
\cite{ra:tran}) which is the identity map on $S$. We also note
that $\phi$ is continuous. The uniqueness of $\phi$ follows from
Theorem~3.4 in \cite{ra:tran}.\hfill $\Box$
\end{proof}

\begin{coro}\label{3.3}$\left.\right.$\vspace{.5pc}

\noindent A triple $(S, T, \circ)${\rm ,} where $(S, T)$ is a
compact Hausdorff space and $(S, \circ)$ a right quasigroup with
identity{\rm ,} can be embedded as a right transversal of a closed
subgroup in a Hausdorff topological group if and only if the map
$\circ$ and the map $\chi$ from $S\times S$ to $S$ defined by
$\chi (x, y)\circ x = y$ are continuous.
\end{coro}

\begin{proof}
Since a compact Hausdorff space is a minimal Tychonoff space, $(S,
T, T', \circ)$ is a topological right transversal if and only if
$T = T'$.\hfill $\Box$
\end{proof}

\begin{coro}\label{3.4}$\left.\right.$\vspace{.5pc}

\noindent $S^{n} = \{x \in \mathbb{R}^{n}{\rm :}\ \Vert x \Vert = 1 \}$
can be embedded {\rm (}topologically{\rm )} as a right transversal
of a closed subgroup in a Hausdorff topological group if and only
if $n = 0, 1, 3, 7$.
\end{coro}

\begin{proof}
Suppose that $S^{n}$ is embedded as a right transversal of a
closed subgroup in a topological group $G$. Then by Corollary~3.3,
there exists a continuous binary operation $\circ$ on $S^{n}$ such
that $(S^n, \circ)$ is a right quasigroup with identity and the
map $\chi\hbox{:}\ S^{n} \times S^{n} \rightarrow S^{n}$ given
by $\chi ((x, y)) \circ x = y$ is continuous. By a well-known
result due to Adams \cite{ad:hinv}, $n = 0, 1, 3$  or $7$.
Conversely $S^{0}, S^{1}$ and $S^{3}$ are Lie groups with respect
to real, complex and quarternionic multiplication respectively.
Let $\circ$ denote the Cayley multiplication on $S^{7}$. Then
$\circ$ is continuous such that $(S^{7}, \circ)$ is a right
quasigroup with identity. Further the map $\chi\hbox{:}\ S^{7}
\times S^{7} \rightarrow S^{7}$ satisfying $\chi ((x, y))
\circ x = y$ is given by $\chi((x, y)) = y \circ x^{-1}$ (note
that Cayley multiplication satisfies $(x \circ y^{-1}) \circ y = x
\circ (y^{-1} \circ y) = x$ and $(y \circ z)^{-1} = z^{-1} \circ
y^{-1}$ although the multiplication is not associative) and so
$\chi$ is also continuous. The result now follows from
Corollary~3.3.\hfill $\Box$
\end{proof}

The following example shows that there exists a proper right
transversal which does not admit any non-discrete locally compact
Hausdorff topology.

\begin{exampl}{\rm \label{3.5}
Let $K$ be an uncountable group which does not admit any
non-discrete $T_1$ topology making $K$ a topological group (for
existence, see p.~339 of \cite{ols:gps} and \cite{she:wp}). Take
$G = K \times D$, where $D$ is a finite group having a non-normal
subgroup $H$ of prime order $p$. Let $T$ be a topology on $G$ with
respect to which it is a $T_1$-topological group. Then as argued
in Example~2.2, $T$ is the discrete topology on $G$. Let $H$ be a
non-normal subgroup of $D$ of order $p$. Then $\{e\}\times H$ is a
non-normal subgroup of $G$ of order $p$. Let $S$ be a transversal
of $\{e\}\times H$ in $G$ which generates $G$. Then $G_{\!S}S \cong
G$. Thus $G_{\!S}S$ admits only discrete topology. Now, if $S$ admits
a non-discrete locally compact Hausdorff topology, then by
Theorem~3.2, we get non-discrete locally compact Hausdorff
topology on $G = G_{\!S}S$.}
\end{exampl}

\section{Transversals of subgroups (not necessarily closed)}

This section is devoted to the following problem.

\setcounter{theore}{0}

\begin{probl}{\rm \label{4.1}
Consider a triple $(S, T,\circ)$, where $(S, T)$ is an infinite
Tychonoff space and $(S,\circ)$ a right quasigroup with identity.
Can we embedd $(S,\circ)$ as a right transversal of a subgroup in
a topological group $G$ so that the subspace topology on $S$ is
$T$?}
\end{probl}

The answer to the above problem is in negative if we stick to
the right transversals of closed subgroups (see Example~4.3). The
following theorem answers the above problem affirmatively,
provided we are satisfied with the right transversal of a subgroup
(not necessarily closed) in a Hausdorff topological group.

\begin{theor}[\!]\label{4.2}
Let $(S, T)$ be a Tychonoff space. Then for any right quasigroup
structure on $S$ with $e$ as the identity element{\rm ,} $S$
becomes a right transversal to a subgroup $H$ {\rm (}not
necessarily closed{\rm )} in a Hausdorff topological group $G$
such that {\rm (i)} the induced right quasigroup structure on $S$
is the given right quasigroup structure and {\rm (ii)} the subspace
topology induced on $S$ in $G$ is the given topology $T$ and $S$
is closed in $G$.
\end{theor}

\begin{proof}
Let $(S, \circ)$ be a right quasigroup with $e$ as the identity.
Let $G = F(S\!\setminus\!\{e \})$ be the free topological group on
the subspace $S\!\setminus\!\{e \}$ of $S$ \cite{mar:ftop}.
Algebraically, $F(S\!\setminus\!\{e\})$ is  simply a free group on
$S\!\setminus\!\{e \}$ \cite{mar:ftop}. We also know \cite{mar:ftop}
that $S\!\setminus\!\{e \}$ is a closed subspace of $G$. Identifying
$1\in G$ with $e \in S$, it is possible to introduce a topology
$T$ on $G$ which induces the given topology on $S$ and in which
$S$ is closed \cite{mar:ftop}. Let $G_{\!S}S$ be the group extension
of the right quasigroup $S$. Since $G_{\!S}S$ is generated by
$S\!\setminus\!\{e \}$, there is a unique surjective homomorphism
$\eta$ (say) from $G$ to $G_{\!S}S$ whose restriction on $S$ is the
identity map on $S$. Thus $G_{\!S}S$ is a topological group with
respect to the quotient topology (induced by $\eta$) which need
not be Hausdorff. The topology on $S$ induced by the topology of
$G_{S}S$ may be strictly coarser than the topology $T$ on $S$ (see
Example~2.3). However, $S$ can also  be considered as a right
transversal of a subgroup $H = \eta^{-1} (G_{S})$ in $G$. If $g
\in G, g \not = 1$, then $\eta (g) = hx$ for some $h \in G_{S}$
and $x \in S$. Since $\eta$ is surjective, there is $k \in H$ such
that $\eta (k) = h$ and then $\eta (g) = hx = \eta (k) \eta(x) =
\eta(kx)$. This shows that $g = k'x$ for some $k' \in H$. Also, if
$g = k_{1}x_{1} = k_{2}x_{2}$, where $k_{1}, k_{2} \in H$ and
$x_{1}, x_{2} \in S$, then $\eta(g) = \eta(k_{1})x_{1} =
\eta(k_{2})x_{2}$. This means $x_{1} =  x_{2}$ and hence $k_{1} =
k_{2}$. One can also observe that the right quasigroup structure
induced on the right transversal $S$ of $H$ in $G$ is the same as
the right quasigroup structure on $S$ from which we started.\hfill
$\Box$
\end{proof}

\begin{exampl}{\rm \label{4.3}
Consider the topological right transversal $S$ discussed in
Example~2.2. Let $T$ be any non-discrete homogeneous Tychonoff
topology. Then as in the notation of the above theorem, $H =
\eta^{-1}(G_S)$ is not closed as otherwise we get a non-discrete
$T_1$-topological group  structure on $G_SS$ in which $S$ becomes
a topological right transversal of the  closed subgroup $G_S$.
This  contradicts Example~2.2.}
\end{exampl}

\section*{Acknowledgement}

We would like to thank the referee for his/her valuable comments
and suggestions.

\end{document}